\documentclass[a4paper]{article}
\usepackage{amsthm,amsmath,amssymb,vmargin}
\usepackage{comment}
\title{\textbf{On the decomposition of the Foulkes module}}                    
\author{Eugenio Giannelli}

\newtheorem{teo1.2}{Theorem 1.2}[section]
\newtheorem{teo1.1}{Theorem 1.1}[section]
\newtheorem{teo}{Theorem}[section]
\newtheorem{cor}[teo]{Corollary}
\newtheorem{lemma}[teo]{Lemma}

\newtheorem{prop}[teo]{Proposition}

\newtheorem*{conjecture}{Conjecture}

\theoremstyle{definition}
\newtheorem{define}[teo]{Definition}

\begin{document}
\bibliographystyle{archmath}

\maketitle

\begin{abstract}
The Foulkes module $H^{(a^b)}$ is the permutation module for the
symmetric group $S_{ab}$ given by the action of $S_{ab}$ on the
collection of set partitions of a set of size $ab$ into $b$ sets
each of size $a$. The main result of this paper is a sufficient 
condition for a simple $\mathbb{C}S_{ab}$-module to have zero multiplicity
in $H^{(a^b)}$. A special case of this result implies that no Specht
module labelled by a hook partition $(ab-r,1^r)$ with $r \ge 1$ 
appears in $H^{(a^b)}$.
\end{abstract}




\section{Introduction}

For $a$ and $b$ natural numbers, let $\Omega^{(a^b)}$ denote the collection of all
set partitions of $\{1,2,\ldots,ab\}$ into $b$ sets each of size $a$. 
Let $H^{(a^b)}$ denote the corresponding $\mathbb{C}S_{ab}$-permutation module, known as the \emph{Foulkes module}.
Let $\phi^{(a^b)}$ be the 
permutation character of $S_{ab}$ afforded by $H^{(a^b)}$.

At the end of Section 1 of \cite{Foulkes}, Foulkes made a conjecture
which can be stated as follows.
\begin{conjecture}[Foulkes' Conjecture]
For all $a$ and $b$ natural numbers such that $a\geq b$,
there exists an injective $\mathbb{C}S_{ab}$-homomorphism from $H^{(a^b)}$ to $H^{(b^a)}$.
\end{conjecture}
The conjecture has been proved to be true only when $b=2$ by Thrall (see \cite{Thrall}), when $b=3$ by Dent (see \cite[Main Theorem]{Dent2000236}), when $b=4$ by McKay (see \cite[Theorem 1.2]{McKay20082050}) and when $b$ is very large compared to $a$ by Brion (see \cite[Corollary 1.3]{springerlink:10.1007/BF03026558}).

Foulkes' original statement of the conjecture was 
as an inequality between multiplicities, namely that,
for all $a$ and $b$ natural numbers such that $a\geq b$ and for all 
partitions  $\lambda$ of $ab$,
$$\left\langle\phi^{(a^b)},\chi^\lambda\right\rangle\leq\left\langle\phi^{(b^a)},\chi^\lambda\right\rangle,$$
where $\chi^\lambda$ is the irreducible character 
of $S_{ab}$ canonically
labelled by $\lambda$.
From this point of view the decomposition of the 
Foulkes module as a direct sum of simple modules becomes central.
 Except in the case when $a=2$ or $b=2$ (see \cite[Chapter 2]{Thrall} and 
\cite{Saxl}) and when $b=3$ (see \cite[Theorem 4.1]{Dent2000236}), 
little is known about the multiplicities of simple modules 
in this decomposition. 
In \cite[Theorem 5.4.34]{JK} an explicit fomula is given for the specific case of simple
 modules labelled by two-row partitions: in this case 
 Foulkes' Conjecture 
 holds with equality. 
We give a short alternative proof of this result in Corollary~\ref{cor:2row} below.
In \cite{WilPag} Paget and Wildon gave a combinatorial description of the minimal partitions that label 
simple modules appearing as summands of Foulkes modules. 

The aim of this paper is to prove a number of new results on when these multiplicities vanish.
We start by giving some standard notation and definitions in Section 2.1. In Section 2.2 we discuss some basic properties of the Foulkes module and we 
describe its restriction to the subgroups $S_r \times S_{ab-r}$ of $S_{ab}$.
In Section 3 we prove the following result which shows that no Specht module labelled by a hook partition $(ab-r,1^r)$ is a direct summand of the Foulkes module $H^{(a^b)}$. 
\begin{teo} \label{hook}
If $a,b$ and $r$ are natural numbers such that $1\leq r<ab$, then  
$$\left< \phi^{(a^b)},\chi^{(ab-r,1^r)} \right> =0.$$ 
\end{teo}

In Section 4 we extend this result, by giving a sufficient condition on a partition $\lambda$ of $ab$ for $\langle\phi^{(a^b)},\chi^\lambda\rangle$ to equal zero.

We need the following notation: let $\alpha=(\alpha_1,\ldots,\alpha_t)$ be a partition of $m\in\mathbb{N}$, let $k\in\mathbb{N}$ 
be 
such that $k\geq t$ and $ab-k-m\geq\alpha_1+1$. Define $[k:\alpha]$ to be the partition 
$$(ab-k-m,\alpha_1+1,\ldots,\alpha_t+1,1^{k-t})$$
of $ab$. (The value of $ab$ will be always clear from the context.)
It is obvious that every partition of $ab$ can be expressed uniquely in the form $[k:\alpha]$. We will call $\alpha$ the \emph{inside-partition} of $[k:\alpha]$.

The main result of this paper is as follows.

\begin{teo}\label{ytt}
Let $a,b$ and $k$ be natural numbers and let $[k:\alpha]$ be a partition of $ab$ with $\alpha=(\alpha_1,\ldots,\alpha_t)$ and $t\leq k$. Let $n:=\sum_{j=2}^{t}\alpha_j$. Suppose that $k > n$ and $\alpha_1 < \frac{1}{2}(k-n)(k-n+1)$. 
Then $$\left<\phi^{(a^b)},\chi^{[k:\alpha]}\right>=0.$$
\end{teo}

Notice that for every simple $\mathbb{C}S_{ab}$-module labelled by $\lambda$, a partition of $ab$ satisfying the hypothesis of Theorem \ref{ytt}, Foulkes' Conjecture holds with equality. Indeed for all $a\geq b$ we have $$\left<\phi^{(a^b)},\chi^{\lambda}\right>=0=\left<\phi^{(b^a)},\chi^{\lambda}\right>,$$
since there is not any restriction on $a$ and $b$ in the statement of the theorem.

By Proposition~\ref{prop:manyparts} below, if
$\langle \phi^{(a^b)},\chi^\lambda \rangle \not= 0$ then $\lambda$ has
at most $b$ parts. When we consider only
characters labelled by such partitions, it occurs that a significant
proportion of the characters appearing with zero multiplicity
in $\phi^{(a^b)}$ satisfy the hypotheses of Theorem~\ref{ytt}.
For example, computations using the computer algebra package
{\sc magma} \cite{Magma} show that there
are $1909$ partitions~$\lambda$ of $30$ with at most $10$ parts 
such that \hbox{$\langle \phi^{(3^{10})},\chi^\lambda \rangle=0$};
of these $492$ satisfy the hypotheses of Theorem~\ref{ytt}.

For an important subclass of partitions to 
which Theorem \ref{ytt} applies we 
refer the reader 
to Corollary \ref{last}.


\section{Preliminaries}

\subsection{Notation and definitions}

A partition $\lambda$ of $n$ is a non-increasing finite sequence of positive integers $$\lambda=(\lambda_1,\lambda_2,\ldots,\lambda_s)$$
such that $\sum_{i=1}^{s}\lambda_i=n$. We write 
$\lambda\vdash n$ to denote that $\lambda$ is a partition of $n$.
The number of parts of a partition $\lambda$ will be denoted by $p(\lambda)$. Denote by $\lambda'$ the conjugate partition of $\lambda$, 
as defined by $\lambda'_j = | \{ i : \lambda_i \ge j \}|$ for $j$ such that $1\leq j\leq p(\lambda)$,
and notice that 
$\lambda_1=p(\lambda')$. We may also denote a partition $\lambda$ by 
$$(\lambda_1^{m_1},\ldots,\lambda_r^{m_r})$$ 
to underline that $\lambda$ has exactly $m_j$ parts equal to 
$\lambda_j$ for all $j\in\{1,2,\ldots,r\}$. 
It will often be useful to think of a partition $\lambda$ as a $\lambda$-Young diagram, as explained by James in \cite[Chapter 3]{james1978representation}. 

\begin{define}
Let $\lambda$ be a partition of $n$ and $\mu$ be a partition of $m$. We say that $\lambda$ is a \emph{subpartition} of $\mu$, and write $\lambda\subseteq\mu$, if $\lambda_j\leq\mu_j$, for all $j$ such that $1 \le j \le \min\bigl( p(\lambda), p(\mu)\bigr)$. 
\end{define}
In particular $\lambda$ is a subpartition of $\mu$ if and only if the Young diagram of $\lambda$ is contained in the Young diagram of $\mu$.
\begin{define}
A \emph{hook partition} is a partition of the form $$\lambda=(n-k,1^k)$$ where $0\leq k\leq n-1$. The number $k$ is called the \emph{leg length} of the hook partition $\lambda$.
\end{define}
We shall also need the \emph{dominance order} on the set of partitions
of a fixed natural number $n$. Given $\lambda, \mu \vdash n$, 
we say $\lambda$ \emph{dominates} $\mu$, and write $\lambda \unrhd \mu$, if  
$$\sum_{i=1}^j\lambda_i\geq\sum_{i=1}^j\mu_i$$
for all $j$ such that $1 \le j \le \min\bigl( p(\lambda), p(\mu) \bigr)$.

Following the definitions and notation of \cite{james1978representation} 
we denote by $M^\lambda$ the Young permutation $\mathbb{C}S_n$-module linearly spanned by the $\lambda$-tabloids and by $S^{\lambda}$ the 
Specht module linearly spanned by the $\lambda$-polytabloids; let $\pi^\lambda$ and $\chi^\lambda$ respectively
denote the associated characters.
From \cite{james1978representation}, we have the following fundamental results.
\begin{teo}\label{james1}
Let $\lambda$ and $\mu$ be two partitions of $n$. If $S^\lambda$ is a direct summand of $M^\mu$ then $\lambda$ dominates~$\mu$.
\end{teo}

\begin{teo}[Branching Theorem]\label{Branch}
Let $\mu$ be a partition of $n$. Let $\Lambda$ be the set of all the partitions of $n+1$ corresponding to the Young diagrams obtained by adding a box to the Young diagram of $\mu$. Then the induced module $S^\mu\big\uparrow^{S_{n+1}}$ decomposes as follows:
$$S^\mu\big\uparrow^{S_{n+1}}=\bigoplus_{\lambda\in\Lambda}S^\lambda.$$
\end{teo} 
The following theorems are straightforward corollaries of 
the Littlewood-Richardson rule, as stated in \cite[Chapter 16]{james1978representation}.
\begin{teo}\label{LR}
Let $k$ be a natural number such that $k<n$ and let $\lambda$ be a partition of $n-k$. If $L$ is the set of all the partitions of $n$ corresponding to the Young diagrams obtained by adding $k$ boxes, no two in the same column, to the Young diagram of $\lambda$, then 
$$(\chi^\lambda\times 1_{S_{k}})\big\uparrow^{S_n}_{S_{n-k}\times S_k}=\sum_{\mu\in L}\chi^\mu.$$
\end{teo}

\begin{teo}\label{LR2}
Let $k$ be a natural number such that $k<n$, let $\lambda$ be a partition of $n-k$, let $\mu$ be a partition of $k$ and let $\nu$ be a partition of $n$. If $$\left\langle(\chi^\lambda\times\chi^\mu)\big\uparrow^{S_n},\chi^\nu\right\rangle\neq 0$$
then $\lambda,\mu\subseteq\nu$, and $p(\nu)\leq p(\lambda)+p(\mu)$.
\end{teo}

\subsection{The Foulkes module}

Here we present some properties of the Foulkes module $H^{(a^b)}$ that will be needed to prove the two main theorems.
  
\begin{prop}\label{prop:manyparts}
Let $\lambda$ be a partition of $ab$ such that $p(\lambda)>b$. Then $$\left\langle\phi^{(a^b)},\chi^{\lambda}\right\rangle=0.$$
\end{prop}

\begin{proof}
It is easily seen that there is a injective map from $H^{(a^b)}$ to $M^{(a^b)}$.
The proposition now follows from Theorem~\ref{james1}.
\end{proof}


\begin{define}
Let $r,a$ and $b$ be natural numbers. We define $P(r)_a^b$ to be the set of all partitions of $r$ with at most $b$ parts and first part of size at most $a$.
\end{define}

An element of $\Omega^{(a^b)}$ can be denoted by $\{A_1,\ldots,A_b\}$, where for each $1\leq j\leq b$, $A_j$ is a subset of $\{1,2,\ldots ab\}$ of size $a$ and for all $i,j$ such that $1\leq i<j\leq b$ it holds $A_i\cap A_j=\emptyset$. 

\begin{define}

Let $r$ be a natural number such that $r<ab$ and let $\lambda$ be in $P(r)_a^b$. We will say that an element $$\{A_1,\ldots,A_b\}\in\Omega^{(a^b)}$$ is \emph{linked} 
to $\lambda$ if the composition of $r$ whose parts are 
$$\bigl|\{1,2,\ldots,r\}\cap A_i\bigr|\quad \text{for $1 \leq i \leq b$}.$$
has underlying partition $\lambda$.
\end{define}
\begin{define}\label{sburt}
Let $r$ be a natural number, such that 
$r<ab$ and let $\lambda$ be in $P(r)_a^b$. We denote by $\mathcal{O}(\lambda)$ the set of all the set partitions in $\Omega^{(a^b)}$ linked to $\lambda$ and by $V^\lambda$ the transitive permutation module for $\mathbb{C}\big(S_{r}\times S_{ab-r}\big)$ linearly spanned by the elements of $\mathcal{O}(\lambda)$.
\end{define}

In the following proposition we show how the restriction $H^{(a^b)} \big\downarrow_{S_{r} \times S_{ab-r}}$ of the Foulkes module decomposes into the direct sum of transitive permutation modules. Such decompositions will be used in all the proofs of the main theorems of this paper.
\begin{prop}\label{QQQ}
Let $r$ be a natural number such that $r < ab$. Then 
\[ H^{(a^b)} \big\downarrow_{S_{r} \times S_{ab-r}} = \bigoplus_{\lambda \in P(r)_a^b} V^\lambda. \]
\end{prop}
\begin{proof}
Let $G = S_{r} \times S_{ab-r}$.
The restriction of $H^{(a^b)}$ to $G$ decomposes as a direct sum of transitive permutation
modules, one for each orbit of $G$ on $\Omega^{(a^b)}$.
Observe that two
set partitions $\mathcal{P}$, $\mathcal{Q} \in \Omega^{(a^b)}$ are in the same orbit
of $G$ on $\Omega^{(a^b)}$ if and only if $\mathcal{P}$ and $\mathcal{Q}$ are linked to the same partition $\lambda \in P(r)_a^b$. The result follows.
\end{proof}

An immediate corollary of Proposition \ref{QQQ} is the following result about the multiplicity of characters labelled by two-row partitions, proved by a different 
argument in \cite[Theorem 5.4.34]{JK}.

\begin{cor}\label{cor:2row}
Let $r,a$ and $b$ be natural numbers. Then 
\begin{itemize} 
\item[(i)] $\bigl\langle\phi^{(a^b)},\pi^{(ab-r,r)}\bigr\rangle=
\bigl\langle\phi^{(b^a)},\pi^{(ab-r,r)}\bigr\rangle=\left|P(r)_a^b\right|$
\item[(ii)] $\bigl\langle\phi^{(a^b)},\chi^{(ab-r,r)}\bigr\rangle=
\bigl\langle\phi^{(b^a)},\chi^{(ab-r,r)}\bigr\rangle=
\left|P(r)_a^b\right|-\left|P(r-1)_{a}^b\right|$
\end{itemize}
\end{cor}
\begin{proof}
It is well known that $$\pi^{(ab-r,r)}=1_{S_{ab-r\times S_{r}}}\big\uparrow^{S_{ab}}.$$
Therefore, by Frobenius reciprocity 
and Proposition \ref{QQQ} we have
\begin{align*}
\left\langle\phi^{(a^b)},\pi^{(ab-r,r)}\right\rangle =&\left\langle\phi^{(a^b)}\big\downarrow_{S_{ab-r}\times S_{r}},1_{S_{ab-r\times S_{r}}}\right\rangle\\
=&\sum_{\lambda\in P(r)_a^b}\left\langle\chi_{V^{\lambda}},1_{S_{ab-r\times S_{r}}}\right\rangle\\
=&\sum_{\lambda\in P(r)_a^b}1=\left|P(r)_a^b\right|.
\end{align*}
To complete the proof of part (i), 
it suffices to observe that the conjugation of partitions induces 
a one to one map between $P(r)_a^b$ and $P(r)_b^a$, and so
$$\left|P(r)_a^b\right|=\left|P(r)_b^a\right|$$ 
for all $r$, $a$ and $b$ natural numbers. 

Part (ii) follows immediately from (i) since $\chi^{(ab-r,r)}=\pi^{(ab-r,r)}-\pi^{(ab-(r-1),r-1)}$.
\end{proof}

We conclude this section with the definition and a description of a generalized Foulkes module that will be used in the proof of Theorem \ref{ytt}.

\begin{define}\label{defn:genFoulkes}

Let $\eta=(a_1^{b_1},\ldots,a_r^{b_r})$ be a partition of $n$, where $a_1>a_2>...>a_r>0$, and let $G = S_{a_1b_1}\times\cdots\times S_{a_rb_r}\leq S_n$. 
We define $$H^\eta=\big(H^{(a_1^{b_1})}\otimes H^{(a_2^{b_2})}\otimes\cdots\otimes H^{(a_r^{b_r})}\big)\big\uparrow_G^{S_n}$$ 
We denote by $\psi^\eta$ the character of the generalized Foulkes module $H^\eta$. 

\end{define}

\begin{define}
Let $\eta=(\eta_1,\ldots,\eta_r)$ be a partition of $n$. Define $\Omega^\eta$ to be the collection of all the set partitions of $\{1,2,\ldots,n\}$ into $r$ sets of sizes $\eta_1,\eta_2,\ldots,\eta_r$.
\end{define} 
The proof of the following proposition is left to the reader.
\begin{prop}\label{general}
Let $\eta=(\eta_1,\ldots,\eta_r)$ be a partition of $n$. Then $H^\eta\cong\mathbb{C}\Omega^\eta.$ 
\end{prop}


\section{The multiplicities of hook characters are zero}
In this section we will prove that no Specht module labelled by a hook partition $(ab-r,1^r)$ appears in the Foulkes module $H^{(a^b)}$.
\begin{define}
Let $U$ be a $\mathbb{C}S_n$-module with character $\chi$. For all $k\in\mathbb{N}$ we denote by $A^kU$ the $k^{th}$ exterior power of $U$, and by $A^k\chi$ the corresponding character.
\end{define}

Let $\epsilon_k$ be the sign character of the symmetric group $S_k$ for any natural number $k$. We leave to the reader the proofs of the following two well known results.

\begin{lemma}\label{sissi}
Let $n$ and $k$ be natural numbers, then $$A^k\pi^{(n-1,1)}=(\epsilon_k\times 1_{n-k})\big\uparrow^{S_n}.$$
\end{lemma}
\begin{lemma}\label{biribissi}
Let  $n$ and $k$ be natural numbers such that $1\leq k\leq n$. Then $$\chi^{(n-k,1^k)}=A^k\chi^{(n-1,1)}\ \ and\ \ A^k\pi^{(n-1,1)}=\chi^{(n-k,1^k)}+\chi^{(n-(k-1),1^{k-1})}.$$
\end{lemma}

In the following proposition we will calculate the inner product between the Foulkes character $\phi^{(a^b)}$ and the character $A^k\pi^{(n-1,1)}$. This is a fundamental step in the proof of Theorem \ref{hook}.
\begin{prop} \label{p1}
Let $a,b$ and $k$ be natural numbers and let $\psi:=\pi^{(ab-1,1)}$. Then 
\label{equ1} 
\[ 
\left<\phi^{(a^b)},A^{k}\psi\right>=
\begin{cases} 0 & \text{if $k\geq 2$} \\
1 & \text{if $k =0, 1$.} 
\end{cases} \]

\end{prop}

\begin{proof}
Firstly consider the case 
$k\geq 2$.
Let $K=S_{\{1,2,\ldots,k\}}\times S_{\{k+1,\ldots,ab\}}\cong S_{k}\times S_{ab-k}\leq S_{ab}$. By Lemma \ref{sissi}
\[
\left<\phi^{(a^b)},A^{k}\psi\right> = 
\left<\phi^{(a^b)},(\epsilon_k\times 1_{ab-k})\uparrow^{S_{ab}}\right> =
\left< \phi^{(a^b)}\!\downarrow_{K},\epsilon_k\times 1_{ab-k}\right>. 
\]
The final inner product above is equal to the number of $\mathbb{C}K$-submodules $U$ in $H^{(a^b)}\downarrow_{K}$ whose associated character is $\epsilon_k\times 1_{S_{n-k}}$.
By Proposition~\ref{QQQ} it suffices to show that if $\lambda\in P(k)_a^b$ then $V^\lambda$ has no submodule with character $\epsilon_k\times 1_{ab-k}$.
Suppose that $u\in V^\lambda$ spans such a submodule. Let $u=\sum_{\mathcal{P}}c_\mathcal{P}\mathcal{P}$, where the sum is over all set partitions $\mathcal{P}\in \mathcal{O}(\lambda)$. Choose $\mathcal{Q}$ such that $c_\mathcal{Q}\neq 0$.

 If $\lambda_1>1$ then there exist $x,y\leq k$ such that $x$ and $y$ appear in the same set in $\mathcal{Q}$. Hence $\mathcal{Q}(x\ y)=\mathcal{Q}$ whereas $u(x\ y)=-u$, a contradiction. Therefore $\lambda=(1^k)$.

If $\lambda=(1^k)$ then   
$$\mathcal{Q}=\bigl\{ \{1,x^1_{2},\ldots ,x^1_{a}\}, 
\{2,x^2_{2}, \ldots ,x^2_{a}\}, \ldots,
\{k,x^k_{2}, \ldots ,x^k_{a}\}, \ldots 
\bigr\}$$
for some $x^i_j \in \{k+1,\ldots, ab \}$. 
Taking $\tau=(1\ 2)(x_1^1\ x_1^2)\cdots (x_a^1\ x_a^2)$, we obtain a contradiction 
again, since $\mathcal{Q}\tau=\mathcal{Q}$ but $u\tau=-u$. 

Hence there are no $\mathbb{C}K$-submodules of $H^{(a^b)}\!\downarrow_K$ 
having character $\epsilon_k\times 1_{S_{ab-k}}$.

The two cases $k=0$ and $k=1$ are left to the reader.
\end{proof}

We are now ready to prove Theorem \ref{hook}. This theorem follows 
at once from Proposition~\ref{p1},
since, from Lemma \ref{biribissi} we have that
$$\chi^{(ab-r,1^r)}=A^r\chi^{(ab-1,1)}=(-1)^r\sum_{k=0}^{r}(-1)^kA^k\psi.$$

We end this section with a corollary of Theorem~\ref{hook} that will be
needed in the proof of Theorem~\ref{ytt}. 
Recall that $\psi^\eta$ is the character of the generalized Foulkes module
$H^\eta$, as defined in Definition~\ref{defn:genFoulkes}.

\begin{cor}\label{Xyzw}
Let $\eta = (a_1^{b_1}, \ldots, a_t^{b_t})$ be a partition
of $n$, where $a_1 > \ldots > a_t$. 
 If $r \ge t$ then
\[ \left< \psi^\eta, \chi^{(n-r,1^r)} \right> = 0. \]
\end{cor}
\begin{proof}
From the definition of generalized Foulkes module, we can write $\psi^{\eta}$ as a character induced from 
\[ \phi^{(a_1^{b_1})} \times \cdots \times \phi^{(a_t^{b_t})}. \]

It follows from Theorems \ref{LR2} and \ref{hook} that
in order  to obtain $\chi^{(n-r,1^r)}$ as an irreducible
constituent of the induced character, we have to take the trivial character in each factor.
Therefore
\[ \left< \psi^\eta, \chi^{(n-r,1^r)} \right> = 
 \left< \big(1_{S_{a_1m_1}} \times \cdots \times 1_{S_{a_tm_t}}\big)\big\uparrow^{S_n}, 
\chi^{(n-r,1^r)} \right>. \]
Observe that
the right-hand side is the multiplicity of $\chi^{(n-r,1^r)}$ in the Young permutation
character $\pi^{(a_1m_1,\ldots,a_tm_t)}$. 
By Theorem~\ref{james1}, 
the constituents of $\pi^{(a_1m_1,\ldots,a_tm_t)}$
are labelled by partitions with at most $t$ parts, 
so we need $t \ge r+1$ to
get a non-zero multiplicity. 
\end{proof}


\section{A sufficient condition for zero multiplicity}
In this section we will prove Theorem \ref{ytt} by an inductive argument. Part of the section will be devoted to the proof of the base step of such induction.
 
Firstly we need to state two technical lemmas.
Let $\beta<ab$ be a natural number. Denote by $K$ the subgroup
$S_{\{1,2,\ldots,\beta\}}\times S_{\{\beta+1,\ldots,ab\}}\cong S_{\beta}\times S_{ab-\beta}$. Let $\lambda$ be in $P(\beta)_a^b$ and let
$V^\lambda$ and $\mathcal{O}(\lambda)$ be as in Definition \ref{sburt}.
Then by a standard result on orbit sums we have the following lemma.

\begin{lemma} \label{didi}
The largest $\mathbb{C}K$-submodule of $V^{\lambda}$ on which $S_{\beta}$ acts trivially is $$U:= \Bigl<\sum_{\sigma\in S_{\beta}}\mathcal{P}\sigma\ |\ \mathcal{P}\in \mathcal{O}(\lambda)
\Bigr>_{\mathbb{C}}.$$
\end{lemma}





With the next lemma we will understand precisely the structure of this particular module $U$.

\begin{lemma} \label{dodo}

Let 
$\lambda\in P(\beta)_a^b$. 
If $\lambda = (\lambda_1, \ldots, \lambda_r)$ then 
$$U\cong \mathbb{C}_{S_\beta}\otimes H^\eta$$
where $\eta=(a^{(b-r)},a-\lambda_r,\ldots,a-\lambda_2,a-\lambda_1)$ and $H^\eta$ is a 
generalized $\mathbb{C}(S_{ab-\beta})$-Foulkes module.

\end{lemma}
\begin{proof}

By Proposition \ref{general} it suffices to show that the set $$\mathcal{W}:=\bigl\{\sum_{\sigma\in S_{\beta}}\mathcal{P}\sigma\ \ |\ \ \mathcal{P}\in \mathcal{O}(\lambda)\bigr\}$$
is isomorphic as a $S_{ab-\beta}$-set to the set $\Omega^{\eta}$ of all $\eta$-partitions of $\{\beta+1,\ldots,ab\}$.

Let $X=\{\beta+1,\beta+2,\ldots,ab\}$. We define a map $f_\lambda:\mathcal{O}(\lambda)\longrightarrow\Omega^\eta$ by $$\mathcal{P}f_\lambda=\{A_1\cap X,A_2\cap X,\ldots,A_b\cap X\}$$ 
where $\mathcal{P}=\{A_1,\ldots,A_b\}$.

It is easy to see that $f_\lambda$ is well defined since 
$\mathcal{O}(\lambda)f_\lambda\subseteq\Omega^\eta$ by definition of $\mathcal{O}(\lambda)$.
The map $f_\lambda$ is surjective, 
and for all $\mathcal{P}$ and $\mathcal{Q}$ in $\mathcal{O}(\lambda)$ we have that $\mathcal{P} f_\lambda=\mathcal{Q} f_\lambda$ if and only if $\mathcal{P}$ and $\mathcal{Q}$ are in the same 
$S_\beta$-orbit of $\mathcal{O}(\lambda)$.
It is easy to see that $f_\lambda$ is an $S_{ab-\beta}$-map 
and that for all $\tau\in S_{\{1,\ldots,\beta\}}$
 we have that $(\mathcal{P}\tau)f_\lambda=\mathcal{P} f_\lambda$, since $\tau$ fixes the numbers greater than $\beta$.

To conclude the proof we define 
$$\tilde{f_\lambda}:\mathcal{W}\longrightarrow\Omega^\eta$$ 
by $$\big(\sum_{\sigma\in S_{\beta}}\mathcal{P}\sigma\big)\tilde{f_\lambda}=\mathcal{P} f_\lambda$$ for all $\mathcal{P}\in \mathcal{O}(\lambda)$.
The map $\tilde{f_\lambda}$ is well defined and 
the surjectivity of $\tilde{f_\lambda}$ follows directly from the surjectivity of $f_\lambda$.
The map $\tilde{f_\lambda}$ is also injective since
$$\big(\sum_{\sigma\in S_{\beta}}\mathcal{P}\sigma\big)\tilde{f_\lambda}=\big(\sum_{\sigma\in S_{\beta}}\mathcal{Q}\sigma\big)\tilde{f_\lambda}\Longleftrightarrow \mathcal{P} f_\lambda=\mathcal{Q} f_\lambda\Longleftrightarrow \mathcal{P}=\mathcal{Q}\tau \Longleftrightarrow\sum_{\sigma\in S_{\beta}}\mathcal{P}\sigma=\sum_{\sigma\in S_{\beta}}\mathcal{Q}\sigma$$
for some $\tau\in S_{\beta}$.

Finally $\tilde{f_\lambda}$ is an $S_{ab-\beta}$-map since $f_\lambda$ is 
an $S_{ab-\beta}$-map  
and 
$\sigma\tau=\tau\sigma$ for all $\sigma\in S_\beta$ and $\tau\in S_{ab-\beta}$.
Therefore $\tilde{f}_\lambda$ is the desired isomorphism.
\end{proof}

In the following proposition we use the notation $[k:\alpha]$ as defined in the introduction. In particular we consider partitions $[k:\alpha]$ of $ab$ with trivial inside-partition $\alpha$ (one row). The proposition is actually the base step of the inductive proof of Theorem \ref{ytt}. 
\begin{prop}\label{Mark}
Let $a,b$ and $k$ be natural numbers. 
For all $\beta<\frac{1}{2}k(k+1)$ we have $$\left\langle\phi^{(a^b)},\chi^{[k:(\beta)]}\right\rangle=0.$$
\end{prop} 
\begin{proof}
By Theorem \ref{LR} and Frobenius reciprocity, we have 
\begin{eqnarray*}
\left\langle\phi^{(a^b)},\chi^{[k:(\beta)]}\right\rangle &\leq &\left\langle\phi^{(a^b)},\big(1_{S_{\beta}}\times\chi^{(ab-(k+\beta),1^{k})}\big)\big\uparrow_{S_{\beta}\times S_{ab-\beta}}^{S_{ab}}\right\rangle\\
&=&\left\langle\phi^{(a^b)}\big\downarrow_{S_{\beta}\times S_{ab-\beta}},1_{S_{\beta}}\times\chi^{(ab-(k+\beta),1^{k})}\right\rangle.
\end{eqnarray*}
Let $K:=S_{\beta}\times S_{ab-\beta}$. By Proposition \ref{QQQ} we have: 

$$H^{(a^b)}\big\downarrow_{K}=\bigoplus_{\lambda\in P(\beta)_a^b}V^{\lambda}.$$



Fix $\lambda=(\lambda_1^{m_1},\lambda_2^{m_2},\ldots,\lambda_{s}^{m_s})\in P(\beta)_a^b$. Let $r:=\sum_{i=1}^s m_i$ be the number of parts of $\lambda$. 
We are now interested in submodules $U\subseteq V^{\lambda}$ such that $S_{\{1,2,\ldots,\beta\}}\cong S_{\beta}$ acts trivially on $U$.
By Lemmas \ref{didi} and \ref{dodo}, 
the largest submodule $U$ of $V^\lambda$ is isomorphic to $H^{\eta}\otimes \mathbb{C}_{S_{\beta}}$, where $\eta=(a^{(b-r)},(a-\lambda_s)^{m_s},\ldots,(a-\lambda_1)^{m_1})$. From now on we will denote $\zeta=((a-\lambda_s)^{m_s},\ldots,(a-\lambda_1)^{m_1})$. Note that 
$$U\cong H^\eta\otimes\mathbb{C}_{S_{\beta}}\cong \big(H^{(a^{b-r})}\otimes H^\zeta\big)\big\uparrow^{S_{ab-\beta}}\otimes \mathbb{C}_{S_{\beta}}.$$
Hence
\begin{align*}
\left\langle\chi_{V^{\lambda}},\chi^{(ab-(k+\beta),1^{k})}\times 1_{S_{\beta}}\right\rangle =&\left\langle\chi_{U},\chi^{(ab-(k+\beta),1^{k})}\times 1_{S_{\beta}}\right\rangle\\
=&\left\langle\big(\phi^{(a^{(b-r)})}\times\psi^\zeta\big)\big\uparrow^{S_{ab-\beta}}\times 1_{S_\beta},\chi^{(ab-(k+\beta),1^{k})}\times 1_{S_{\beta}}\right\rangle\\
=&\left\langle\big(\phi^{(a^{(b-r)})}\times\psi^\zeta\big)\big\uparrow^{S_{ab-\beta}},\chi^{(ab-(k+\beta),1^{k})}\right\rangle\\
=&\sum_{\nu,\mu}d_{\nu}^{\mu}\left\langle(\chi^\nu\times\chi^\mu)\big\uparrow^{S_{ab-\beta}},\chi^{(ab-(k+\beta),1^{k})}\right\rangle 
\end{align*}
where $\chi^{\nu}$ is an irreducible character of $S_{a(b-r)}$ with non zero multiplicity in $\phi^{(a^{(b-r)})}$, 
$\chi^{\mu}$ is an irreducible character of $S_{ar-\beta}$ having non zero multiplicity in $\psi^\zeta$, and $d_{\nu}^\mu$ is the multiplicity of their tensor product in the decomposition of $H^\eta$. 
Notice that the last sum is not equal to zero if and only if there exist $\nu$ and $\mu$ such that $(\chi^\nu\times\chi^\mu)\big\uparrow^{S_{ab-\beta}}$ 
 contains a $hook$ character of $S_{ab-\beta}$ having leg length equal to $k$ in its decomposition. By Theorem~\ref{LR2}, we have that both $\nu$ and $\mu$ must be subpartitions of $(ab-(k+\beta),1^{k})$. 
 This means that $\nu$ and $\mu$ are hooks or trivial partitions. In particular we deduce from Theorem \ref{hook} that $\nu=(a(b-r))$. So we need $\mu$ to be a hook with leg length at least $k-1$ to have $$\left\langle(\chi^\nu\times\chi^\mu)\big\uparrow^{S_{ab-\beta}},\chi^{(ab-(k+\beta),1^{k})}\right\rangle\neq 0.$$
 On the other hand $$\psi^\zeta=\big(\phi^{((a-\lambda_1)^{m_1})}\times\cdots\times\phi^{((a-\lambda_s)^{m_{s}})}\big)\big\uparrow^{S_{ar-\beta}}$$
So by Corollary \ref{Xyzw} we have that the hooks that have non-zero multiplicity in the decomposition of $\psi^\zeta$ have at most $s$ parts, where $s$ is the number of different parts of $\lambda$. 

We observe that the smallest number $\tilde\beta$ having a partition $\lambda$ with $k$ different parts is $\frac{k(k+1)}{2}$, with $\lambda=(k,k-1,\ldots,2,1)$.
So under our hypothesis $\beta<\frac{k(k+1)}{2}$ we obtain that $\chi^{\mu}$ cannot be a hook character with leg length at least $k-1$.
Hence for all $\lambda\in P(\beta)_a^b$ we have that $$\left\langle\chi_{V^{\lambda}},\chi^{(ab-(k+\beta),1^{k})}\times 1_{S_{\beta}}\right\rangle=0.$$
\end{proof}

We are now ready to prove Theorem \ref{ytt}.

\begin{proof}[Proof of Theorem~\ref{ytt}]
We proceed by induction on $t$, the number of parts of the inside-partition $\alpha$. 

If $t=1$ then $$\left\langle\phi^{(a^b)},\chi^{[k:(\alpha_1)]}\right\rangle=0$$ by Proposition \ref{Mark}.

Suppose now that $t>1$ and the theorem holds when the inside-partition has less then $t$ parts. Denote $\nu:=(\alpha_1,\alpha_2,\alpha_3,\ldots,\alpha_{t-1})$.
By Theorem~\ref{LR}, Lemma \ref{dodo} and Frobenius reciprocity 
we have that 
\begin{align*}
\left\langle\phi^{(a^b)},\chi^{[k:\alpha]}\right\rangle 
\leq& \left\langle\phi^{(a^b)},\big(\chi^{[k:\nu]}\times 1_{S_{\alpha_t}}\big)\big\uparrow^{S_{ab}}\right\rangle\\
=&\left\langle\phi^{(a^b)}\big\downarrow_{S_{ab-\alpha_t}\times S_{\alpha_t}},\chi^{[k:\nu]}\times 1_{S_{\alpha_t}}\right\rangle\\
=&\sum_{\lambda\in P(\alpha_t)_a^b}\left\langle\chi_{V^\lambda},\chi^{[k:\nu]}\times 1_{S_{\alpha_t}}\right\rangle\\
=&\sum_{\lambda\in P(\alpha_t)_a^b}\left\langle\chi_{U^\lambda},\chi^{[k:\nu]}\times 1_{S_{\alpha_t}}\right\rangle\\
=&\sum_{\lambda\in P(\alpha_t)_a^b}\left\langle\big(\phi^{(a^{(b-p(\lambda))})}\times\psi^{(a-\lambda_{p(\lambda)},\ldots,a-\lambda_{1})}\big)\big\uparrow^{S_{ab-\alpha_t}},\chi^{[k:\nu]}\right\rangle\\
=&\sum_{\lambda\in P(\alpha_t)_a^b}\big(\sum_{\zeta,\mu}d^\lambda_{\zeta\mu}\left\langle\big(\chi^\zeta\times\chi^\mu\big)\big\uparrow^{S_{ab-\alpha_t}},\chi^{[k:\nu]}\right\rangle\big),
\end{align*}
where, for each $\lambda\in P(\alpha_t)_a^b$, $U^\lambda$ is the largest $\mathbb{C}(S_{ab-\alpha_t}\times S_{\alpha_t})$ submodule of $V^\lambda$
on which $S_{\alpha_t}$ acts trivially and $\sum_{\zeta,\mu}d_{\zeta\mu}^\lambda(\chi^{\zeta}\times\chi^{\mu})$ is the decomposition into irreducible characters of the character
$\phi^{(a^{(b-p(\lambda))})}\times\psi^{(a-\lambda_{p(\lambda)},\ldots,a-\lambda_{1})}$. 

For every $\lambda\in P(\alpha_t)_a^b$, 
observe that every simple summand $S^\mu$ of $H^{(a-\lambda_{p(\lambda)},\ldots,a-\lambda_{1})}$ is a simple summand of the Young permutation module $M^{(a-\lambda_{p(\lambda)},\ldots,a-\lambda_{1})}$. Hence by Theorem \ref{james1} we have that the partition $\mu$ has at most $p(\lambda)$ parts; in particular it has at most $\alpha_t$ parts. It follows that, by Theorem~\ref{LR2}, we need $\zeta$ to have at least $k+1-\alpha_t$ parts, and to be a subpartition of $[k:\nu]$ in order to have 
$$\left\langle\big(\chi^\zeta\times\chi^\mu\big)\big\uparrow^{S_{ab-\alpha_t}},\chi^{[k:\nu]}\right\rangle\neq 0.$$
Therefore $\zeta$ must be of the form $$[k_{\zeta}:\beta]\vdash a(b-p(\lambda))$$
with   
\begin{itemize}
\item $\beta=(\beta_1,\ldots,\beta_{s})\subseteq \nu$, and
\item $k_{\zeta}\geq k-\alpha_t$.
\end{itemize} 
We conclude proving that such a $\zeta$ cannot label any irreducible summand of the Foulkes character $\phi^{(a^{(b-p(\lambda))})}$.
 
Define $n_\zeta:=\sum_{j=2}^{s}\beta_j$; if $s=1$ then let $n_\zeta = 0$.
We observe that such a partition $\zeta$ has inside-partition $\beta$ having $s\leq t-1$ parts and it satisfies the initial hypothesis, since
\begin{itemize}
\item $k_{\zeta}\geq k-\alpha_t> n-\alpha_t\geq n_\zeta$, and
\item $\beta_1\leq\alpha_1<\frac{(k-n)(k-n+1)}{2}\leq\frac{(k_{\zeta}-\sum_{j=2}^{t-1}\alpha_j)(k_{\zeta}-\sum_{j=2}^{t-1}\alpha_j+1)}{2}\leq\frac{(k_{\zeta}-n_\zeta)(k_{\zeta}-n_\zeta+1)}{2}$.
\end{itemize}
Hence $\chi^\zeta$ has zero multiplicity in $\phi^{(a^{(b-p(\lambda))})}$ by induction. Therefore 
\begin{align*}
\left\langle\phi^{(a^b)},\chi^{[k:\alpha]}\right\rangle 
\leq&\sum_{\lambda\in P(\alpha_t)_a^b}\big(\sum_{\zeta,\mu}d^\lambda_{\zeta\mu}\left\langle\big(\chi^\zeta\times\chi^\mu\big)\big\uparrow^{S_{ab-\alpha_t}},\chi^{[k:\nu]}\right\rangle\big)=0.
\end{align*}
The theorem is then proved.
\end{proof}

As mentioned in the introduction, and as we will prove in the following corollary, a consequence of our main theorem is 
that every Specht module labelled by a partition having leg length equal to~$k$ and at most $k$ boxes inside the hook has zero multiplicity, except when the $k$ boxes are \emph{column-shaped} (i.e. the inside-partition is $(1^k)$). In that particular case we are able to prove that the multiplicity equals 1, for all the values of $k< b$. The proof is similar to that of Proposition \ref{Mark} and is omitted. In \cite[Lemma 3.3]{Dent2000236} Dent proved the same result
in the specific case $k=b-1$.

\begin{cor} \label{last}
Let $a,b,k$ and $m$ be natural numbers. Let $m\leq k$ and $\alpha$ be a partition of $m$ not equal to $(1^k)$. Then
$$\left\langle\phi^{(a^b)},\chi^{[k:\alpha]}\right\rangle=0$$
\end{cor}
\begin{proof}
Let $\alpha$ be an arbitrary partition of $m$ not equal to $(1^k)$. Then $$\alpha=(\alpha_1,\alpha_2,\ldots,\alpha_t).$$
Write $n:=\sum_{j=2}^t\alpha_j$. We will show that $\chi^{[k:\alpha]}$ satisfies the hypothesis of Theorem \ref{ytt} that:
\begin{itemize}
\item $k>n$, and
\item $\alpha_1<\frac{1}{2}(k-n)(k-n+1)$. 
\end{itemize}
The first condition is trivial since $$k\geq m=\alpha_1+n.$$
To prove the second condition proceed by contradiction: suppose that $$\alpha_1\geq\frac{(k-n)(k-n+1)}{2}.$$
Then $$k-n\geq\frac{(k-n)(k-n+1)}{2}.$$
This implies $k-n=0$ or $k-n=1$.
The first situation is impossible because  
$0=k-n\geq\alpha_1>0$.
The second is also impossible because 
$0<\alpha_1\leq k-n=1$ implies $\alpha_1=1$ and $\alpha_1+n=k$ with $\alpha=(1^k)$. 
\end{proof}



\providecommand{\bysame}{\leavevmode\hbox to3em{\hrulefill}\thinspace}
\providecommand{\MR}{\relax\ifhmode\unskip\space\fi MR }
\providecommand{\MRhref}[2]{%
  \href{http://www.ams.org/mathscinet-getitem?mr=#1}{#2}
}
\providecommand{\href}[2]{#2}

\end{document}